\documentclass{article}[12pt]
\usepackage{color,times,amsmath,amsfonts,latexsym,epsfig,epsf,graphics,colordvi}
\usepackage{url,hyperref}
\usepackage[english]{babel}

\newcommand{\RR}{\mathbb {R}}

\newcommand{\al}{\alpha }
\newcommand{\bt}{\beta }

\newcommand{\la}{\lambda }

\newcommand{\bp}{\begin{pmat}}
\newcommand{\ep}{\end{pmat}}
\def\tcr{\textcolor[rgb]{1,0,0}}
\def\tcb{\textcolor[rgb]{0,0,1}}

\makeatletter
\def\adots{\mathinner{\mkern1mu\raise\p@\vbox{\kern7\p@\hbox{.}}\mkern2mu\raise4\p@\hbox{.}\mkern2mu\raise7\p@\hbox{.}\mkern1mu}}
\makeatother
 
 \title{\vspace*{-2mm} The Challenges of Teaching\\
 Elementary Linear Algebra in\\  a Modern Matrix Based Way}
 
\begin{document}

\date{}
\author{}
\maketitle

\begin{center}
\vspace*{-11mm}
{\LARGE Frank Uhlig}\\[1mm]
{\large \emph{Department of Mathematics and Statistics, \\ Auburn University, Auburn, AL 36849-5310, USA\\[1mm] \quad {\bf{ uhligfd@auburn.edu}}}}
\end{center}

\vspace*{-1mm}
\begin{center} { \bf Abstract  } \\[2mm]
\begin{minipage}{110mm}

We assess the situation of our elementary Linear Algebra classes in the US holistically and through personal history recollections. Possible remedies for our elementary Linear Algebra's teaching problems are discussed and a change from abstract algebraic taught classes to a concrete matrix based first course is considered. The challenges of such modernization attempts for this course are laid out in light of  our increased after-Covid use of e-books and e-primers.\\[0.9mm]
 We specifically address the useless and needless, but ubiquitous use of determinants, characteristic polynomials  and polynomial root finding methods that are propagated in our elementary text books and are used in the majority of our elementary Linear Algebra classes  for the matrix eigenvalue problem but that have no practical use whatsoever and offer no solution for finding matrix eigenvalues. \\[0.9mm]
 This paper challenges  all mathematicians  as we  have misinformed and miseducated  our students badly for decades in elementary Linear Algebra now and urges a switch to a new, fully matrix theoretical approach that covers all classical subjects in a practical and computable way. 
 \end{minipage}\\[5mm]
\end{center}
 
 The paper is an amalgan of talks given over the last few years on using matrices to teach modern Linear Algebra. It extracts and reformulates many ideas that were first voiced by the author in math education zoom talks  in England and in the US during Covid times. At that time my  main  aim was to exemplify several  modern matrix based approaches to standard  Linear Algebra concepts concretely. These zoom talks  included first trial runs of the current Lesson Plans 3, 5 and 6 that are now part of the ILAS Education Webpages, see [3]. More recent attempts at modernizing our elementary Linear Algebra courses were likewise presented at ILAS Galway and in JMM in Boston in 2022, this time  in person.  Our main aim here is to wrap up the theory of why we need to change elementary Linear Algebra classes completely and how this can be done. We  need to learn   how to effect this change quickly, rightfully and amiably among all the involved actors, namely the instructors, the students, our client disciplines, textbook printing houses and so forth. Such challenges are unprecedented in mathematics and science since Galileo's time some 400 years ago. But they are absolutely necessary  now for Linear Algebra after 200 years of our own stagnation in 1820s to 1870s  mode. \\ 
 How can we, I and you  help this transition?

 \section* {Introduction}

This paper approaches our teaching efforts for elementary beginning Linear Algebra classes around the globe in a holistic manner, looking at this complex problem from a personal history viewpoint and as a community effort.\\[-2mm]

{\em \bf A Personal History and Insights Gained}\\[-3mm]

For many decades, 'Linear Algebra" has been an undergraduate  cornerstone class taken early by  students in mathematics  and those in many applied fields of science, in engineering and other disciplines. I have taught Linear Algebra at the university level for over 50 years now and I have witnessed many pedagogical advances on how to improve the student learning and  subject understanding.
Looking back to the 1960s when I entered the University of Cologne in Germany to study math it was an honor for the chair of a math institute to teach the incoming first year students and the Linear Algebra lectures were given from hand written notes and proofs of the professor's best understanding and interpretation of the subject.\\
 Consequently, there were many individual attempts to build a Linear Algebra course from matrix theory\linebreak
  {[} Emil Artin (Hamburg Lecture Notes, 1961); Prof.  Wolfram  Jehne, Dr. Plewe (Lectures, U K\"oln, 1965); Hans Schneider and Phil Barker, (\emph{Matrices and Linear Algebra}, Dover, 1968) {]} then that were competing with more classical abstract  and theoretical deductive and  logic based Definition-Theorem-Proof and determinant based lectures. In fact, when we as students searched the math library at the Mathematics Institute in Cologne University  for Linear Algebra books, we found most of our answers in books titled as 'Theory of Determinants' or slight variations thereof,  where all questions, proofs and results for nonsingular matrices $A$, for example, were couched in $\det(A) \neq 0$ statements. Searching the card catalogue for "Linear Algebra"  or walking along the 50 yards of math books, 6 feet high, and searching by eye for 'Linear Algebra' on their spines brought no results then and we were stuck with interpreting determinantal results.\\[1mm]
In the 1970s I learnt and mostly taught  numerical analysis and returned to teach elementary Linear Algebra by 1980. At that time I noticed disturbing ripples and a discontent in the student body when they tried to grasp  abstract math and algebraic proofs that were the vogue then in top-down teaching Definition-Theorem-Proof style classes, both in Germany and in the US.\\[-2mm]

{\em \bf The Development of Linear Algebra Teaching Insights and the Coursing of Elementary Linear Algebra Classes}\\[-3mm]

At the same time I began to become aware of  efforts to adapt our pedagogy to better understand our students' learning difficulties (Jean Piaget, Anna Sierpinska, Guershon Harel, Tommy Dreyfuss and others). Linear Algebra often appears as a multicolored chameleon with identical symbols used to mark different items such as the number zero (0) or the zero vector (0) or a trivial subspace, described by $\{0\}$. Such notational vaguenesses lead to confusion and ultimately to a  lack of interest and engagement in Linear Algebra students who are taught the subject by rote and tote row reductions, determinant evaluations, characteristic polynomial computations and so forth, but never given the key to the kingdom of understanding the subject deeply.  Our textbooks do not describe the ways, needs and advances of Linear Algebra and Matrix Theory as our subject is being used more and more in modern applications and has become a foundation for the internet, in modern commerce, engineering, and innovations. In the US, Linear Algebra was mostly then and still is taught as a service course by young faculty, post docs and other faculty who - in the times of  'publish or perish' - have to develop their own specific fields of expertise to gain tenure and promotions and not by experts who know and would promote the modern ways with Linear Algebra. Educational studies of Linear Algebra subjects eventually led to  more sensible ways to reach and teach our students and to small changes in our syllabi and textbooks. But the student malaise has survived and after two Linear Algebra Education Study Groups (1993 and 2022) not much has improved. Except that dedicated teachers now teach more and more individually from applications and on-line through experimentations with matrices rather than top-down Definition-Theorem-Proof wise.\\[1mm]
 Teaching is a complex process. 
 Piece-meal answers to complex problems generally result in unintended consequences. Therefore a holistic analysis of any complex problem that we encounter when teaching is necessary and much preferred. The holistic approach to a  complex problem starts with a simple, somewhat personal question for the concerned:
How would you want the situation in question to be resolved and become irrelevant in the future, such as in 5, 10 or 30 years? We want to let elementary Linear Algebra courses run smoothly, with the current problems completely removed. For  Linear Algebra  this  means: a complete understandings of the  class'   objective, of the what and hows, the whys; and completely up to date in subject matter  taught and balanced out between all classroom actors. A holistic solution of the elementary Linear Algebra teaching problem would remove the discontent of our students and replace those ill feelings with natural curiosity and personal discovery. This would  replace our current disturbing class presentations that are so at odds  with our modern  understandings and daily applications of matrices.\\[2mm]
Some of us might want to bring software into the first Linear Algebra course, others on-line course-ware and so forth. Before talking of specific remedies, we will now look at the holistic net that ties teachers and students together.\\[-6mm]

\section{The Web of Forces and Realms that Govern the Mathematical Teaching Process and the First Linear Algebra Course}

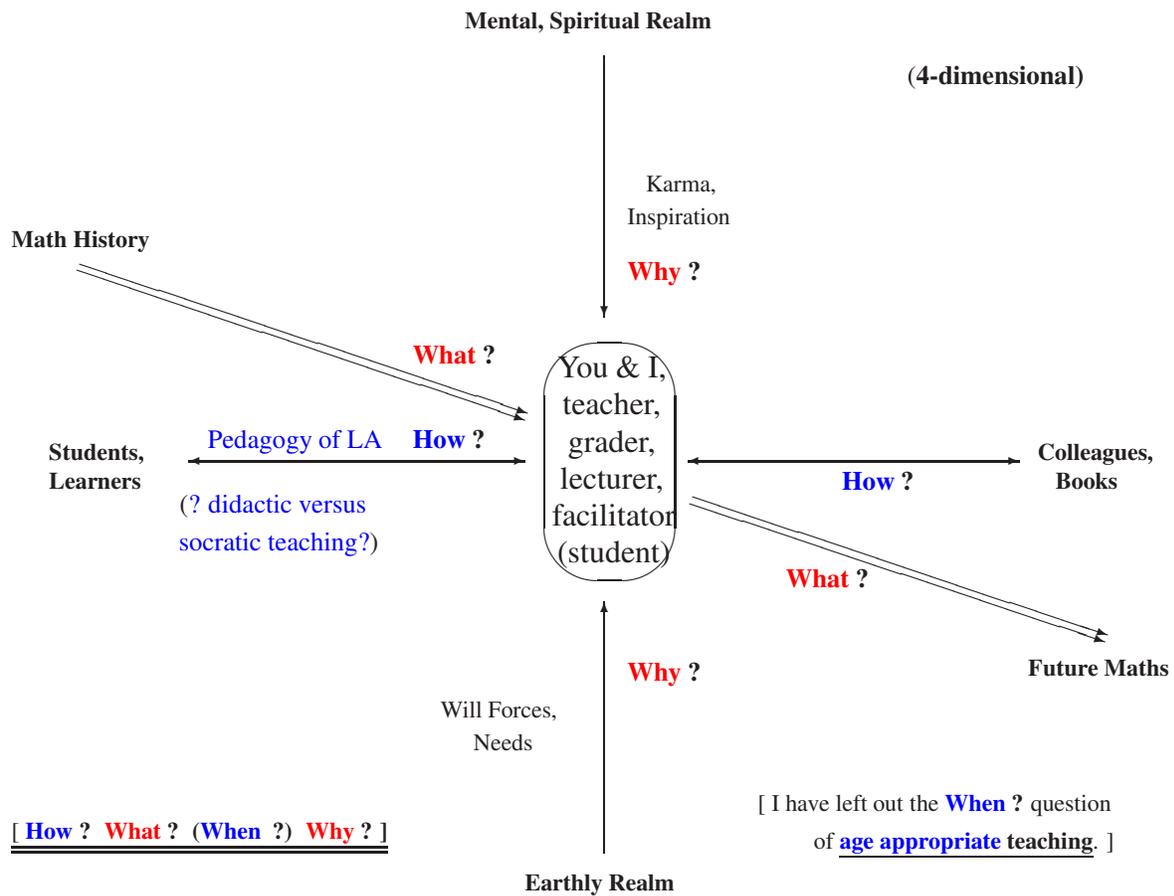
\setlength{\unitlength}{1.24mm}

\vspace*{2mm}
\hspace*{-2mm}\begin{picture}(128,98)(0,0)


\put(60.5,49.8){\small \begin{minipage}{14.5mm} {\large \begin{center} You \& I,\\teacher,\\grader,\\lecturer,\\$\!$facilitator\\ (student)\end{center} }\end{minipage}}

\put(6,50){\small \begin{minipage}{11mm} \begin{center} \bf Students,\\ Learners \end{center} \end{minipage}}
\put(22,50.5){\vector(1,0){35}} \put(36,50.5){\vector(-1,0){15}}
\put(45,61){\bf \tcr{What} ?}
\put(23,52){\tcb {Pedagogy of LA}}  
\put(20,45){(\tcb{? didactic versus}}   
\put(16,41){\tcb{ \hspace*{4mm} socratic teaching?})} 

\put(45,52){\bf \tcb{How} ?} 

\put(91,47.5){\bf \tcb{How} ?}
\put(85,37){\bf \tcr{What} ?}  

\put(112,50){\small \begin{minipage}{13mm} \begin{center} \bf Colleagues,\\ Books\end{center} \end{minipage}}

\put(75,50.5){\vector(1,0){35}} \put(89,50.5){\vector(-1,0){14.5}}

\put(50.6,97){\small \bf Mental, Spiritual Realm}
\put(70,79.5){\small Karma,}
\put(68,76){\small Inspiration}
\put(68,70){\bf \tcr{Why} ?}
\put (65.5,94){\vector(0,-1){28}}
\put(98,91){(\bf 4-dimensional)}   

\put(57,4.5){\small \bf Earthly Realm}
\put(48,23){\small Will Forces,}
\put(51.5,19.5){\small Needs}
\put(68,27){\bf \tcr{Why} ?}
\put (65.5,8.5){\vector(0,1){27}}

\put(2,73.5){\small \bf{Math History}}
\put(9,71){\vector(3,-1){48}}
\put(9.3,71.7){\vector(3,-1){48}}

\put(111,27.5){\small \bf{Future Maths}}
\put(74.7,46){\vector(3,-1){44.5}}
\put(75,46.7){\vector(3,-1){44.5}}

\put(82,13){\small [ I have left out the  \tcb{\bf When} {\bf ?} question}
\put(88,9){\small   of {\bf \underline{\tcb{age appropriate} teaching}}. ]}

\put(2,10){\small  \underline{\underline{[ \bf \tcb{How} ? \  \tcr{What} ?  \ (\tcb{When } ?) \ \tcr{Why} ? ]}} }

\put(66.1,50.5){\oval(14,25.5)}
\end{picture}\\[-8mm]

\setlength{\unitlength}{1mm}

The most fundamental question of any personal holistic teaching analysis is 'why do I teach, why do we teach,  why do you teach'. Each teacher has to be conscious of his or her answer to this question every day. What motivates us to teach? There are earthly reasons and karmic ones. We all need to have food and shelter In the modern world we learn and study, then hire ourselves out to gain our 'daily bread'.\\
 Do I, do we or you teach because we must or  is 'teaching' just our job for earning our 'bread'.  Would I or you rather frame houses all day long  or repair cars, or garden  or perform other fruitful work that really excites me, you or us? \\[2mm] 
Besides practical reasons, there are karmic ones such as 'I teach because I am meant to teach', or 'I am called to teach the young', or 'Teaching is my life's work' and so forth because 'teaching makes me happy'.\\[1mm]
 On the daily level there are certain demands on 'How' we teach in a College or University environment. When we start our careers, we may like to use textbooks written by others  or there are prescribed  syllabi that  we might have to follow  as the University or Department has decreed. There are pedagogic principles to follow, as well as teaching styles such as interactive, open classrooms or top-down lecturing and multiple ways of inverse teaching in between.\\[1mm]
When we teach a living, a lively expanding subject such as Linear Algebra, there is its mathematical history from antiquity to the present and into the future that binds us  into  additional  dependencies. Do we teach from antiquity on or start with modern, even cutting-edge new methods, with current applications and thoughts or rather present the old? The same holds on the students' side: what have they been prepared for; at age 8 and in eighth grade? When is a college student's  math comprehension well founded, when are the concepts of Linear Algebra overbearing and how so. These questions connect teachers and students to a daily 'How' to proceed quandary.  Our personal voyages into mathematical thinking  - both as teacher and student - pull us, intersect in us, bind us in our studies. They surface and resonate in us  as actors of the teaching and learning 'game' of College education.\\[2mm]
Here are two lists of linear algebra subjects, associated names and dates from eternities past until near today.\\[2mm]

\hspace*{-2mm}\begin{tabular}{lll}
Gaussian elimination, Gaxby and Saxby & Assur and Babylon& 4000 to 2000 BE\\[1.5mm]

polynomial &Descartes & 1637 \\
polynomial roots &&\\
complex numbers&Cardano&1545\\
Rule of Signs& Descartes& 1637\\[1.5mm]

matrix&Leibniz, Seki&1683\\
determinant&Cardano, Seki, Leibniz&1545, 1683, 1693\\
Cramer's rule &Cramer&1750\\[1.5mm]

Fundamental Theorem of Algebra& Gauss& 1799\\
Cauchy-Schwarz inequality&Cauchy, Schwarz&1821, 1888\\
normal equation& Gauss & 1822\\
Ruth-Hurwitz&Ruth, Hurwitz&1876, 1895\\
Gram-Schmidt&Gram, Schmidt&1883, 1907\\[1.5mm]

eigenvector&Euler&1751, 1760\\
eigenvalue&Cauchy&1819\\
characteristic polynomial&Cauchy&1819\\
Cayley-Hamilton Theorem&Cayley, Hamilton&1858, 1864\\
Jordan normal form&Jordan&1870\\[3.5mm]
\end{tabular}

\vspace*{-85mm}
\hspace*{123mm}
\rotatebox{90}{$\underbrace{\hspace*{80mm}}$}\\

\vspace*{-62mm}\hspace*{129mm}\begin{minipage}[t]{28mm} Standard math contents with algebraic proofs or pencil and paper explanations in 95\% + of our Linear Algebra textbooks and classes today
\end{minipage}\\

\newpage

Of special importance for a holistic approach to teaching Linear Algebra is the diagonal Math History line in the holistic diagram of forces. What we teach in a specific math course that is  fundamental  for many client disciplines determines how relevant the class becomes for our students and how deeply they want to engage in our  classes. Teachers have to be aware of what methods are being used today in  applications. The evolving topics of modern Linear Algebra and matrices below govern how  modern our approach can and should be. Are we preparing our students for a future or just teaching its past? Compare the two lists above and below.\\[1mm]

\hspace*{-2mm}\begin{tabular}{lll}
minimal polynomial& Krylov& 1931\\
vector iteration& Krylov& 1931\\[1.5mm]

LR factorization&&\\
QR factorization&&\\
Schur decomposition&Schur&1909\\
Hessenberg matrix&Hessenberg&1943\\
Givens reflection& Givens&1950s\\
Householder transform&Householder&1958\\
QR algorithm& Francis, Kublanovskaya& 1961\\
Singular Value Decomposition&Golub, Kahan&1965\\
SVD and 'PageRank'  search engine&Page, Brin&1996\\
QR multi-shift algorithm&Braman, Byers, Mathias&2002\\
TensorFlow search engine&Google team& 2021\\[1.5mm]
\end{tabular}\\

\vspace*{-58.5mm}\hspace*{103mm}
\rotatebox{90}{$\underbrace{\hspace*{55mm}}$}\\

\vspace*{-39.3mm}\hspace*{110mm}\begin{minipage}[t]{36mm} Math items occasionally mentioned in 5\% $\pm$ of our textbooks or classes.
\end{minipage}\\[25mm]

On the previous page  there is a list of named old and ancient  Linear Algebra subjects that appear in our current elementary Linear Algebra textbooks and are generally taught in our current first elementary  Linear Algebra course. Its 18 braced items  determine the majority of our classroom syllabi of today together with Linear Algebra's standard vector space and specific matrix studies and their inter-relations. Any linear algebraic advances after the 1930s  in the second modern subject list  above are generally not mentioned, let alone taught in the majority of our current elementary Linear Algebra classes.\\[-2mm]

 Unfortunately, most applied  programs in the US do not suggest or require advanced or numerical Linear Algebra courses for their graduate and undergraduate students. This omission  then reverberates through much of the US engineering world. Tech companies must reeducate our graduates upon their first employment to become up to date in Linear Algebra. This is needless and very costly. We must now act  to modernize our first Linear Algebra course and transmit  modern matrix advances of the last 90, 40 or 20 years adequately to our young.\\[0.5mm]
 Cauchy's then brilliant 1820s idea to use determinants and the characteristic polynomial $f(x)$ as eigendata holders for $A_{n,n}$ has never born fruit; humankind was never able to find useful ways to extract eigenvalues from the determinant based characteristic polynomial $f(x) = \det (A - x I_n) = 0$ in 200 years of trying. Why do we continue to teach it?\\[1mm]
 In 1931 Aleksei Krylov  published his seminal paper on vector iteration and Krylov subspaces for matrices.\\[1mm]
{\em From Krylov's paper, from Russian at https ://mathshistory.st-andrews.ac.uk/Biographies/Krylov\_Aleksei/  : \\[1mm]
"It is clear that, if for $k = 2$ and $k = 3$ it is easy to compose this [secular] equation} (i.e., the characteristic polynomial $f(x)$ of a $k$ by $k$ matrix), {\em then for $k = 4$ the laying-out becomes cumbersome, and for values $k$ more than 5 this is completely unrealizable in a direct way.\\[0.5mm]
  Therefore one should use methods where the full development of the determinant is avoided. The aim of the paper ... is to present simple methods of composition of the secular equation in the developed form, after which, its solution, i.e. numerical computation of its roots, does not present any difficulty".}\\[1mm]
Krylov's first paragraph above gives a realistic description of trying to compute the characteristic equation of even very low dimensional matrices $A_{k,k}$. Our current students or 95 \% or more of them become aware of this today in elementary Linear Algebra courses and textbooks.  However, they are taught no eigenvalue developments beyond Cauchy's time. Rather than using determinants to try to find matrix eigenvalues, Krylov  suggested to find a linear dependency of the first $k+1$ vector iterates $\{ b, Ab, A^2b, ..., A^kb\}$ for $b \neq o_n$ and then he describes the minimal polynomial of $A$ that has all of $A$'s eigenvalues as roots. Krylov's latter hope of finding polynomial roots numerically has - however -  never been realized in 200 years of trying.\\[1mm]
Yet  the 'highlight' of our elementary Linear Algebra courses today is Cauchy's dead-end or cul-de-sac idea. This is unbelievable. We have had excellent numerical matrix eigenvalue finders  in Francis' QR algorithm since the 1960s and for matrix dimensions $k$ into the 10,000s with Francis' multi-shift algorithm since the early 2000s, but we and our textbooks most often do not even teach orthogonal matrix factorizations such as QR.\\[-5mm]
 
\section{Mathematics for a Syllabus Change by Using Modern Matrix Theory in Elementary Linear Algebra  Courses Today}

 I have always been attracted to constructive  proofs and specifically to matrix based insights for  linear algebra. In my decades as a  teacher and research mathematician, I have kept aware of certain jewels of matrix based thinking. In 2020/21 I was approached by Ian Benson who had read my earlier textbook titled 'Transform Linear Algebra', see [1]. I had chosen  this provocative title (with 'Transform' either a verb or a noun) to express the power of linear transforms, aka matrices; both  in Linear Algebra itself and when teaching Linear Algebra. Linear (in-)dependence of a set of vectors is  easily accessible from a row echelon form reduction of the associated column vector matrix that make our ubiquitous textbook verbal logic definitions and treatment of linear (in-)dependence laughable. Ian wanted to help British teachers modernize the a-levels university exams by allowing  Linear Algebra as a substitute subject for the Calculus requirement, working through  through the Association of Teachers of Mathematics (ATM). When I looked at  the a-level test primers from Oxford and Cambridge Universities I was amazed that their preparatory Calculus examples and problems were identical to what I had studied and was tested on for my Abitur in Germany in the mid 1960s. Upon asking my father I learned that he had been tested on exactly the same calculus examples in order to pass his Abitur  in 1919/20. Had nothing changed in Mathematics since then? What about modern Matrix Theory?\\[1mm] 
 And I had arrived at a mathematics challenge: Can matrices  help us to understand Linear Algebra better and revolutionize our Linear Algebra syllabi to become relevant in our internet, smartphone and computer times? \\[2mm]
These thoughts propelled me to develop a set of matrix based  Lesson Plans for a modern Linear Algebra course, leading  from millennias old row reduction to eigenvector/eigenvalue computations and beyond that is based on modern practices and  notions  in our second subject, name and date list  from the 1930s on.  \\[2mm]
 Linear Algebra is the study of linear functions $f$ that map $n-$space to $m-$space.
A function $f : u \in \RR^n \to f(u) \in \RR^m$ is called linear if 
$$ f(\al u + \bt v) = \al f(u) + \bt f(v) $$
for all vectors $u$ and $v$ and all scalars $\al$ and $\bt$. Using Riesz Representation Theorem from 1907, we show in our first Lesson Plan that every linear transform $f : \RR^n \to \RR^m$ can be described by its standard matrix representation 
$$A_{m,n} =  \bp \vdots && \vdots \\ f(e_1) & \cdots & f(e_n) \\ \vdots && \vdots \ep $$ 
for the standard unit vectors $e_i \in \RR^n$ that contain zeros in every position, except for a 1 in position $i$. Thus studying linear transformations in Lesson Plan 1 is synonymous with studying matrices.\\[1mm]
In Lesson Plan 2 we study linear equations and matrix row echelon form reductions, thereby discovering pivots and free columns or free variables. Then we interpret these results and determine criteria for solvable and unsolvable systems of linear equations $Ax = b$. Efficient and correct row reduction must be performed by hand, pencil and paper first. Once that procedure is mastered, students can turn to software and learn to code row reduction correctly in Python, in Octave, Mathematica or Matlab and others.\\[1mm] 
Lecture Plan 3 first deals theoretically with matrix inversion, then uses  paper and pencil computations for practical matrix inversion and finally via student built software codes. These techniques are applied to vector spaces, sets of vectors, spanning sets, and the concepts of subspaces, bases and dimensions. Linear independence is casually touched here, too.\\[1mm]
Our linear algebra Lesson Plans have begun with the standard unit vector basis ${\cal E}  = \{e_i\}$ and the standard matrix representation $A_{\cal E}$ of a linear transform  $f$. In Lesson Plan 4 we generalize these concepts to arbitrary bases ${\cal U} = \{u_i\}$ of $n-$space and $A_{\cal U}$, the ${\cal U}$ basis representation of $A_{\cal E}$.
We aim to find a basis  ${\cal U}$ for $A$ that contains the eigenvectors $u_i \neq o_n$ of $A_{n,n}$, i.e., a basis $\{u_i\}$ of $n-$space with $A u_i = \la_i u_i$ for $ i =1, ..., n$. \\[1mm]
Eigenspaces and matrix eigenvalues are then computed in Lesson Plan 5 by using Krylov vector iteration. We also prove that every square matrix does have eigenvalues and eigenvectors by giving an existence proof for the minimal polynomial for $A$ and using the Fundamental Theorem of Algebra.\\[1mm]
Our Lesson Plans end with a study of angles and orthogonality in space, with planar rotations, the trigonometric identities for sine and cosine and orthogonal matrices in Plan 6.\\[1mm] Lesson Plan 7 revisits Gaussian elimination in the form of  an LR factorization of matrices and then expands  matrix factorizations to orthogonal QR decompositions via Householder transformations and simple applications to solve unsolvable systems of linear equation in a least squares sense.\\[1mm]
This set of Lesson Plans gives  one  coherent mathematical path to teach an elementary Linear Algebra class subjects - and a little bit more - completely from and through modern Matrix Theory. Each Lesson Plan takes between 2 and 6 class days.\\
Other such paths (see e.g. Gilbert Strang, \emph{Introduction to Linear Algebra}, 6th edition)   exist that may also be able to free Linear Algebra courses from Cauchy's 200 years old untenable ideas.

\section{Deeper Challenges to Changing our Linear Algebra Syllabi; how can this be done}

The challenge of  mathematically  restructuring our elementary Linear Algebra offerings has been met by this set of Lesson Plans that treat Linear Algebra through Matrix Theory.\\
 But the open challenge to actually transform our entry level Linear Algebra  courses in the US on the human, the personal level seems almost insurmountably more difficult and complex. With humans having a tendency to cling to decades or centuries old ways in  thinking, doing  and feeling, it is very difficult and takes time to enter and establish ourselves consciously in a new computational epoch., see [2] for example. Our new epoch of individuality has been all around us for several decades now, but it has not fully taken power of our spirits. And this is a far deeper challenge than the purely mathematical one was.\\[-2mm]

How easy can  and will it be for faculty to recognize the changes that we are  called  to make now and switch to teach modern methods in our elementary  Linear Algebra classes. How long will it take us to clearly dismiss Cauchy's two centuries old dismal approach as fruitless? When will we abandon it. How soon will other areas such as engineering recognize that the Gram-Schmidt orthogonalizing process - for example - is unstable and ill-advised, but still demand this  'knowledge' of their students in exams.
\\[2mm]
Will students switch on their own and study on-line e-books that use matrix methods instead of the outdated subjects that  might still be taught? Can graders and instructors accept  student work that might prefer and actually use Lesson Plan type matrix knowledge and matrix methods to solve  'classic' abstract Linear Algebra or matrix eigenvalue problems? Should primers be created and deposited on-line that show how to solve standard Linear Algebra problems and still understand or 'speak' the old redundant  language and solve the old test questions via Lesson Plan e-book style matrix methods. Can a client based,  a grass-roots student driven revitalization of our elementary Linear Algebra curriculum really happen? In real time, quickly and thoroughly; I hope so and I know it can.\\[2mm]
Will there be a grass-roots movement somehow that solves 'classical' abstract linear independence problems for example  concretely via row echelon form reductions and gives us quick concrete answers to abstractly formulated  problems in our centuries old  textbooks via software and thereby un-stifles our classrooms?\\[2mm]
Who will switch first, teachers or students? Locally, regionally, or globally. How can we justify deliberately teaching outdated material today? Just  as we have been taught unjustifiably by our teachers who repeated their teachers' inadequate understandings and back and back for a century plus. Can we  finally break this misbegotten circle, this dereliction of our academic  duty to teach to the best of our subject matter's knowledge and  change our teachings to every new generation's knowledge base? The QR algorithm for matrix eigen computations is 60 years old now. Current Linear Algebra research forges ahead with tensor studies and parameter-varying matrix problems for robots, control, with AI, and so forth. And this and much more will have to be explained in elementary Linear Algebra classes in the future.\\[-2mm]

Given the power of matrix methods, there is room, the desire and also the need for a second elementary Linear Algebra course for students of mathematics, in engineering  and all science majors now. Our current students need to be taught, learn and  know so much more about matrices to succeed in their professions.

\section{How to Help the Transition of Teachers and Students from Abstract Linear Algebra to Linear Algebra and Matrices}

\vspace*{2mm}
Modern Lesson Plans that emphasize Matrix Theory as a tool to understand  Linear Algebra  are not enough by themselves. How will a student or teacher switch when all around him or her the talk is about Cauchy's determinant idea and characteristic polynomial root finding techniques, of Cramer's rule, Gram-Schmidt and the like -- from more than a century ago and utterly fruitless in modern applications. There is need for Matrix theoretic explanations of the standard linear algebraic notions and classical 'solution' paths. An alternate  'Lesson Plan* set' of sorts is needed to mitigate and translate  between the two ways of thinking and which gives practical  translations from one viewpoint to the other. A first try is below:

\subsection{Linear Independence and more}

\vspace*{2mm}
Here we want to translate basic linear algebra concepts and proofs that are usually given as problems in abstract Linear Algebra classes and try to solve them from a Matrix viewpoint.\\[1mm]

{\bf (a) } {\em Standard abstract formulation}: Given $k$ vectors $z_i$ in $\RR^n$. There exists one vector $z_\ell$  that is a linear combination of the vectors $\{z_i | i \neq \ell\}$ if and only if the set of vectors $z_i$ is linearly dependent.\\
A classical course's proof  expresses $z_\ell = \sum_{i \neq \ell} \al_i z_i$ and rewrites this equation as $\sum_{i \neq \ell} \al_i z_i -  1 \cdot z_\ell = 0 \in \RR^n$ which makes the $z_i$ linearly dependent by definition because the coefficient of $z_\ell$ is nonzero while the linear combination of the $z_i$ is the zero vector.\\[1mm]
{\em Matrix theoretic approach}: The engine for all  linear (in)dependence questions here is the row echelon $R_{n,k}$ form of the matrix $A_{n,k}$ with columns $z_i$. The invariants of any REF $R$ of $A$ are the row and column locations of its pivots and the location of its free columns when we do not allow column swaps. \\
If the $z_i$ are linearly dependent then the REF $R$ of the column vector matrix $A_{n,k}$ for the $z_i$ has at least one free column in one position $\ell \leq k$ and the linear system  $A_{n,\ell-1} = z_\ell$ can be solved, putting $z_\ell \in \text{ span}\{ z_1, ..., z_{\ell -1}\}$. \
This problem was easy.\\

{\bf (b) }  {\em A more difficult problem in three parts:}\\
 If a set of $m$ vectors $\{z_1,...,z_m\}$ spans a linear subspace $\cal V$ of $\RR^n$ and $z_i \in \text{ span}\{z_j | j \neq i\}$, then the set of $m-1$ vectors $\{z_1 + z_i, ... , z_{i-1} + z_i,z_{i+1} + z_i, ...,z_m+z_i\}$ span $\cal V$.\\
1) \ {\em True} or {\em False}? Proof  or counter example ?\\
2) \ What can be concluded when the given set of $m$ vectors $\{z_1,...,z_m\}$ contains only distinct vectors?\\
3) \ What can happen if the set of $m-1$ vectors $\{z_1 + z_i, ... , z_{i-1} + z_i,z_{i+1} + z_i, ...,z_m+z_i\}$ then becomes\\ \hspace*{4mm} a basis for $V$?\\

{\bf (c) } {\em The multi-roles of {'dimension'} :} \\
1) \ Our colloquial use of 'dimension' refers to the size of a table, of a bed, a door or of a package.\\
2) \ For vectors we often call the number of their entries their dimension. For example a real vector $x$ has\linebreak
\hspace*{4mm} dimension $n$ if it contains $n$ real number entries.\\
3) \ Matrices $A$ have row and column dimensions, referring to the number of rows and columns.\\
4) \ A vector space $V$ has dimension $k$ if any of its bases contains $k$ vectors.\\

{\bf (d) } {\em The multi-role of the symbol '0', called {\em 'zero' :}}\\
List at least six distinct uses of~ $0$ in mathematics.\\
Repeat and try to find distinct uses for 1 or 'one' in math.\\

{\bf (e) } {\em The size of a linear independent spanning set :} \\
{\em Standard abstract formulation}:  If $\ell$ vectors $ w_1, ..., w_\ell$ span a linear subspace $\cal W$ of $\RR^n$ and $S$ is a set  of $m > \ell$ vectors in $\cal W$ then the $m$ vectors in $S$ are linearly dependent.\\[1mm]
{\em Matrix theoretic approach:} The row echelon form of the column vector matrix for the $\ell$ vectors $w_i$ contains exactly $k = \text{ dim}({\cal W})$ pivots and obviously $k \leq \ell$. The REF of the column vector matrix for the set of $m$  vectors in  $S$  cannot have more than $k = \text{ dim}({\cal W})$ pivots. But  $m > \ell \geq \text{ dim}({\cal W}) = k$, so there must be linearly dependent vectors among the $m$ vectors in $S$.\\

{\bf (f) } {\em Linearly independent and spanning vectors in a linear subspace $\cal U$ of $\RR^n$:}\\[1mm]
{\em Standard abstract formulation}:  For a linear subspace ${\cal U} \subset \RR^m$ of dimension $k$, any linear independent  set of $k$ vectors in ${\cal U}$ spans ${\cal U}$ and conversely any set of $k$ vectors in ${\cal U}$ that spans the subspace ${\cal U}$ of dimension $k$ is linearly independent.\\[1mm]
{\em Matrix theoretic approach:}  We again look at the REF of the column vector matrix $A_{m,k}$  for $k = \text{ dim}({\cal U})$. \\[0.5mm]
1) \ We start from $k$  linearly independent vectors $\{ u_1, ..., u_k\}$ in ${\cal U}$ : Can these vectors span the whole subspace $\cal U$? I.e., can every linear system $A_{m,k}x_k = b_m $   for the column vector matrix $A$ of the $u_j$ and $b = b_m \in {\cal U}$ be solved?

As the $k$ columns $u_j$ of $A$ are assumed to be linearly independent, the REF of $A$ has $k$ pivots since no column is a linear combination of the others, see part (a). If there are $k+1$  pivots in the augmented matrix $(A|b)$ then $b \notin {\cal U}$ because the linear system $Ax = b$ is unsolvable. Otherwise if $b \in {\cal U}$ the system is solvable, i.e., for any $ b \in {\cal U}$, the right hand side vector  $b$  is a linear combination of the vectors $\{ u_1, ..., u_k\}$. And therefore the set of vectors $\{u_i\}$ spans $\cal U$.\\[0.5mm]
2) \ If  $k$  vectors $ u_1, ..., u_k$ span the $k$-dimensional subspace $\cal U$ then the REF of their column vector matrix $A_{m,k}$ must have $k$ pivots since otherwise the augmented matrix $(A|b)$ would indicate non-solvability of $Ax=b$ for some $b \in {\cal U}$ and thus their failing to span $\cal U$. With $k$ pivots in $A$ there can be no linear dependency among the $u_i$, see (a). \\[1mm]
{\em Note a certain duality here}: linearly independent sets of vectors in a subspace $\cal V$ of dimension $k>0$ may contain at most $k$ vectors, while spanning sets of vectors for a subspace of dimension $k>1$ must contain at least $k$ vectors. When these vector containment numbers  coincide at $k$, then both the set of $k$ spanning vectors and the set of $k$ linearly independent vectors form a basis for this subspace $\cal V$ of dimension $k$.\\

{\bf (g) } {\em The determinant of a square matrix:}\\[1mm]
Historically the determinant of a square matrix $A_{n,n}$ was defined abstractly in the 16th and 17th century as the {\em sum of signed products} of $n$ entries in mutually exclusive row and column positions of $A$.
$$
\det(A) = \sum _{\sigma \in \Sigma_n} \text{ sign}(\sigma) a_{1,\sigma(1)} \cdot a_{2,\sigma(2)} \cdot  ... a_{n,\sigma(n)} 
$$
where $\Sigma_n$ is the group of all permutations of the integers from 1 to $n$. $\Sigma_n$ contains $n! = 1 \cdot 2 \cdot 3 \cdot ... \cdot n$ permutations. This definition of the determinant has been simplified by {\em row or column expansions} to evaluate determinants and certain rules of the behavior of determinants have been established over time, such as\\[1mm]
\hspace*{15mm} $ \det(A\cdot B) = \det(A) \cdot \det(B)$ \hspace*{7.4mm} for any two $n$ by $n$ matrices $A$ and $B$,\\[0.5mm]
\hspace*{15mm} $\det(E_{j,k}) = -1$ \hspace*{27.7mm} for exchanging rows $j$ and $k$ in $A$,\\[0.5mm]
\hspace*{15mm}  $\det(E_{\al_k,k} \cdot A)  = \al_k~ \det(A)$ \hspace*{9.1mm} for multiplying the row $k$ of $A$ by $\al_k$, \\[0.5mm]
\hspace*{15mm} $\det(\det(E_{(\al_k,k) +j} \cdot A)  = \det(A)$ \quad for adding $\al_k$ times row $k$ to row $j$ of $A$, and\\[0.5mm]
\hspace*{15mm} $\det(T) = \prod t_{i,i}$ \hspace{26.5mm} for any upper or lower triangular matrix $T$ with diagonal\linebreak
\hspace*{67.8mm} entries $t_{i,i}$.\\[1mm]
These five classical determinant rules allow us to use the LR factorization of $A_{n,n} $ to evaluate $$\det(A) = \det(L~R) = (-1)^{ex} ~\det(L) ~ \det(R) = (-1)^{ex}~\prod r(i,i)$$
 with the diagonal entries $r_{i,i}$ of $R$ and $ex = \# \text{ of row interchanges in the LR factorization process of A} $. Here $\det(L) = 1$ because $L$'s diagonal entries are all equal to 1.\\[1mm]
 The determinant of a square matrix $A_{n,n}$ in fact measures the volume of the parallelepiped bordered by the column vectors of $A$ as edges of a squished 'cigarette box' in $n$-space. If the columns of $A$ are linearly dependent, this 'squished box' lies in a flat, less than $n$-dimensional subspace of $\RR^n$ and subsequently has $n$-volume zero. \\[0.5mm]
To try to evaluate $f(x) = \det (A- x I_n)$ as Cauchy set out to do 200 years ago and find the eigenvalues of $A_{n,n}$ where the characteristic polynomial $f$ has its roots never succeeded beyond $n = 4$ or beyond, except for specialized doctored matrices.\\[4mm]

{\large \bf \em Therefore we need to drop Cauchy's unsuccessful determinantal and polynomial roots approach of centuries ago.\\[2mm] 
And  instead adapt (and adopt) modern Matrix Theory into our elementary Linear Algebra courses and follow Krylov and Francis' paths of vector iteration or matrix factorization methods to learn how to extract matrix eigendata, quickly and  accurately from any given matrix; in our computer, engineering,  and software age today.}\\[3mm]
{\Large And teach what we now know.}\\[4mm]

{\Large \bf  Some References}\\[3mm]
[1] \ Frank Uhlig, Transform Linear Algebra, Prentice-Hall, 2002,  ISBN 
0-13-041535-9, 502 + xx p.\\[2mm]
[2] \ Frank Uhlig, The Eight Epochs of Math as Regards Past and Future Matrix Computations, in \emph{Recent Trends in Computational Science and Engineering}, S. Celebi (Ed.),  InTechOpen (2018),\\  {\url  {http://dx.doi.org/10.5772/intechopen.73329}}, 25 p. \\
 {[} complete with graphs and references at \ {\url  {http://arxiv.org/abs/2008.01900}} \ (2020), 19 p. ]\\[2mm]
 [3] Frank Uhlig and Rachel Quinlan, A Library of Lesson Plans for Teaching Linear Algebra, ILAS Image, (69) Fall 2022, p. 3-4.  \\  {\url  {https://31g343.p3cdn1.secureserver.net/wp-content/uploads/image69.pdf}}

\vspace*{22mm}
\hspace*{3mm} in  \ \ \  ...Box/local/latex/MadridILAS23Edu/PrimusEdu23 paper.tex\\[6mm]
 \hspace*{8mm} {\em \today}

\end{document}